\documentclass[11pt]{article} 
\usepackage{amsmath, amsthm, amsfonts, amssymb, latexsym}

\theoremstyle{plain}

\newtheorem{theorem}{Theorem}
\newtheorem{proposition}{Proposition}
\newtheorem{lemma}{Lemma}
\newtheorem{definition}{Definition}
\newtheorem{corollary}[theorem]{Corollary}

\begin{document}
\title{\textbf{Transversely Hessian foliations and information geometry}}
\author {Michel Nguiffo Boyom\\{\small I3M UMR CNRS 5149
Universit\'e Montpellier 2, Pl E. Bataillon F34095 Montpellier CEDEX5}
\and Robert A. Wolak \\{\small Wydzial Mathematyki i Informatyki, Uniwersytet Jagiellonski}\\{\small Lojasiewicza 6,  Krakow, Poland, robert.wolak@im.uj.edu.pl.}}

\maketitle
\begin{abstract} A family of probability distributions parametrized by an open domain $\Lambda$  in $R^n$ defines the Fisher information matrix on this domain which  is positive semi-definite. In  information geometry the standard assumption has been that  the Fisher information matrix tensor   is positive definite defining in this way a Riemannian metric on $\Lambda$. If we replace the "positive definite" assumption by the existence of a suitable torsion-free connection, a foliation with a transversely Hessian structure appears naturally. In the paper we develop the study of transversely Hessian foliations in view of applications in information geometry.

\vskip 2mm
{\footnotesize\noindent{\em$2011$ Mathematics Subject Classification\/}.  53C12 }

\end{abstract}
\vskip 2mm




\maketitle


\section{Introduction --  Information Geometry}
The Fisher metric is one of the basic tools in  information geometry. It is defined on an open domain in $R^m$ which parametrizes the set of probability distributions under consideration. The assumpsions made lead to the study of affine manifolds with Riemannian metrics, with both geometrical structures loosely related. Using the probability distributions one can define a $(0,2)$ tensor field on the open domain. The assumption that the defined tensor field is positive definite is rather strong, therefore we propose to study the consequences of a weaker condition.  In this situation a foliation with a very particular transverse structure appears. Such foliated manifolds are natural generalizations of Hessian manifolds. 

\medskip

Let $\Lambda$ be a domain in ${ R}^m.$ We consider  families of probability distributions on a set $\mathcal X$ parametrized by $\lambda \in \Lambda .$

$$ {\mathcal P} = \{ p(x;\lambda ) \vert \lambda \in \Lambda \}$$

(1) $\Lambda$ is a domain in ${ R}^m,$

(2) $p(x; \lambda ) $ for a fixed $x$ is a smooth function in $\lambda ,$

(3) the operation of integration with respect to $x$ and differentiation with respect to $\lambda$ are commutative.

\medskip

\noindent
{\bf Definition}
{\it Let $ {\mathcal P} = \{ p(x;\lambda ) \vert \lambda \in \Lambda \}$ be a family of probability distributions on a set $\mathcal X$ parametrized by $\lambda \in \Lambda .$ We set  $l_{\lambda} = l(x;\lambda ) = log p(x; \lambda )$ and denote by $E_{\lambda}$ the expectation with respect to $p_{\lambda} = p(x; \lambda).$ Then the matrix $g_F (\lambda ) = [g_{ij}( \lambda )]$ defined by

$$ g_{ij}( \lambda ) = E_{\lambda }[\frac{\partial l_{\lambda}}{\partial \lambda ^i}\frac{\partial l_{\lambda}}{\partial \lambda ^j}]\\
= \int_{\mathcal X}\frac{\partial l(x;\lambda)}{\partial \lambda ^i}\frac{\partial l(x;\lambda)}{\partial \lambda ^j}p(x;\lambda ) dx $$

\noindent
is called the Fisher information matrix tensor. }

\medskip

Simple calculations show, see \cite{Shima}, that

$$g_{ij}( \lambda ) = - E_{\lambda }[\frac{\partial^2 l_{\lambda}}{\partial \lambda ^i \partial \lambda ^j}].$$

\noindent
The Fisher information matrix tensor  $g_F (\lambda ) = [g_{ij}( \lambda ) ]$ is positive semi-definite on $\Lambda$:

$$ \Sigma_{i.j} g_{ij}(\lambda )c^ic^j  = 
\int_{\cal X} \{ \Sigma_i c^i \frac{\partial l(x;\lambda)}{\partial \lambda^i} \}^2 p(x,\lambda )dx \geq 0.$$

\noindent
In  information geometry the standard assumption has been, cf. \cite{Shima}, p.105,:

\medskip

(4)  For a family of probability distributions $ {\mathcal P} = \{ p(x;\lambda ) \vert \lambda \in \Lambda \}$ the Fisher information matrix tensor $g_F (\lambda )= [g_{ij}( \lambda ) ]$ is {\bf positive definite } on $\Lambda .$

\medskip

\noindent
We  {\bf weaken this condition } assuming only that the Fisher information matrix tensor  is a tensor field parallel with respect to  some torsion-free connection on $\Lambda.$ Then a foliation appears in a very natural way, and under   some mild assumptions it has a transverse  Hessian structure. The main part of this note is devoted to the development of the foundations  of the theory of transversely Hessian foliation which can be applied to a classification of spaces of probability distributions in the non-regular case.

\section{Foliations}

 Let $\mathcal F$ be a foliation on an $m$-manifold $M.$ Then $\mathcal F$ is defined by a cocycle
   $\mathcal U =\{ U_i , f_i , k_{ij} \}_{i\in I}$ modeled on a $q$-manifold $N_0$ ($0<q<m$) such that

(1) $\{ U_i \}_{i\in I}$ is an open covering of M, 

(2) $f_i : U_i \to N_0$ are submersions with connected fibres, 

(3) $k_{ij} : N_0\to N_0 $ are local diffeomorphisms of $N_0$ with $f_i = k_{ij} f_j$ on $U_i \cap U_j .$ 

\noindent The connected components of the trace of any leaf of $\mathcal F$ on $U_i\;$ consists of fibres of $f_i .$ The open subsets $N_i = f_i ( U_i )\subset N_0$ form a $q$-dimensional manifold $N_{\mathcal U}=\bigsqcup N_i$, which can be considered to be  a complete transverse manifold of the foliation $\mathcal F .$ The pseudogroup $\mathcal H_{\mathcal U}$ of local diffeomorphisms of $N$ generated by $k_{ij}$ is called the holonomy pseudogroup of the foliated manifold $(M, \mathcal F )$ defined by the cocycle $\; \mathcal U .$  The equivalence class ${\mathcal H}$ of ${\mathcal H}_{\mathcal U},$  for the notion of the pseudogroup equivalence see  \cite{Ha1,Ha2,Ha3}, is called the holonomy group of ${\mathcal F},$ or of the foliated manifold $(M,\mathcal F )$. A foliation on a smooth manifold $M$ understood as an involutive subbundle of $TM, $ or equivalently, according to the Frobenius theorem, cf. \cite{CC}, p.37, as a partition of the manifold by submanifolds of the same dimension with some regularity condition, can be defined by many different cocycles. There is a notion of equivalent cocycle, modeled on the notion of equivalent atlases of a smooth manifold, and a foliation can be understood as an equivalence class of such cocycles, cf., the notion of a smooth structure on a topological space. Moreover, a pseudogroup equivalent to a holonomy pseudogroup representative is itself a pseudogroup asssociated to some cocycle defining the foliation.   Therefore in some cases, for our foliation, we will be able to choose  a cocycle modeled on a particular  manifold, cf. \cite{Wo0}.

\medskip

The vector bundle $N(M,{\mathcal F})=TM/T{\mathcal F}$ is called the normal bundle of the foliation $\mathcal F$. Then the tangent bundle $TM$ is isomorphic to the direct sum $T{\mathcal F} \oplus N(M,{\mathcal F})$. These isomorphisms are determined by the choice of a supplementary subbundle $Q$ in $TM$ to the tangent bundle to the foliation $T{\mathcal F}$. The cocycle   $\mathcal U =\{ U_i , f_i , k_{ij} \}_{i\in I}$ modeled on a $q$-manifold $N_0$  induces on the normal bundle a cocycle   $\mathcal V =\{ V_i ,\bar{ f}_i , \bar{k}_{ij} \}_{i\in I}$ modeled on the $2q$-manifold $TN_0,$ where $V_i=TU_i,$ $\bar{f}_i$ is the mapping induced by $df_i,$ and $\bar{k}_{ij} =dk_{ij}$. The foliation ${\mathcal F}_N$  of the normal bundle is of codimension $2q,$ its leaves project on leaves of $\mathcal F$. They are, in fact, coverings of these leaves. In a similar way one can foliate any bundle obtained via a point-wise process from the normal bundle, e.g., the frame bundle of the normal bundle, the dual normal  bundle, any tensor product of these bundles.

\medskip

  Let  $\phi : U\to R^p \times R^q ,\,\phi =(\phi^1 ,\phi^2)=(x^1, \cdots , x^p ,y^1,\cdots ,y^q )$ be an adapted chart on a    foliated manifold $(M, \mathcal F)$. Then on $U$ the vector fields
   $\frac{\partial}{\partial x^1},\cdots, \frac{\partial}{\partial x^p}$ span    the bundle $T\mathcal F$ tangent to the leaves of the foliation $\mathcal F ,$ the equivalence
   classes of $\frac{\partial}{\partial y^1}, \cdots ,   \frac{\partial}{\partial y^q}$ denoted by $\bar{\frac{\partial}{\partial y^1}}, \cdots, \bar
   {\frac{\partial}{\partial y^q}}$, span the normal bundle $N(M, \mathcal F ) $,  in fact, these vector fields are foliated sections of the normal bundle foliated by the foliation ${\mathcal F}_N.$

\medskip
 
In the case of a foliated manifold we can consider three types of geometrical structures related to the foliation:

\medskip

\noindent
{\it transverse} - defined on the transverse manifold, the associated holonomy pseudogroup consists of automorphisms of this geometrical structure;

\medskip
\noindent
{\it foliated} -  only defined on the normal bundle, and when expressed in a  local adapted chart, depending only on the transverse coordinates; a foliated structure projects to a transverse structure along submersions of the cocycle defining the foliation;

\medskip
\noindent
{\it associated} - defined globally, on the tangent bundle but adapted to the spliting, and defining a foliated structure on the normal bundle.

\medskip
Foliated and transverse structures are in one-to-one correspondence, an associated structure defines a foliated structure, but different associated structures can define the same foliated structure, cf. \cite{Wo0}.

\medskip


\medskip
\noindent
{\bf Example 1} Let us see what  these types of structures  give in the case of a Riemannian metric on a foliated manifold $(M,{\mathcal F})$ with the foliation defined by a cocycle ${\mathcal U}$. A transverse Riemannian metric is a Riemannian metric $\hat{g}$ on the transverse manifold $N$ of which elements of the holonomy pseudogroup are local isometries. If it is true for one holonomy group representative, it is true for any  equivalent pseudogroup. Such a Riemannian metric can be "lifted" to a metric tensor field $g$ on the normal bundle. Locally, using the sections  $\bar{\frac{\partial}{\partial y^1}}, \cdots, \bar    {\frac{\partial}{\partial y^q}}$ and their duals $dy_1,...,dy_q$ such a tensor field can be written as

$$g= \Sigma_{i=1}^q g_{ij}(y)dy_idy_j$$

\noindent
Mind that the condition "foliated" is equivalent to the fact that the functions $g_{ij}$  depend only on the variables $y^1,...,y^q.$
Of course, a metric tensor on the normal bundle can be extended to a metric tensor on the tangent bundle $TM$ using the splitting and the isomorphism we have discussed above. The normal bundle is isomorphic to any complementary subbundle $Q$ of the tangent bundle $T{\mathcal F}$. Therefore the tensor field $g$ can be transported to a metric tensor field $g_Q$ on the subbundle $Q$. Making $Q$ orthogonal to $T{\mathcal F}$ and choosing a metric tensor  in $T{\mathcal F}$ we get a Riemannian metric $\bar{g}$ on $M,$ which induces a metric tensor $g$ in the normal bundle. Such Riemannian metrics are called {\it bundle-like\/} and  have very interesting properties and characterisations, cf. \cite{Molino}.

\medskip
\noindent
{\bf Example 2} Likewise there three types of connections "adapted" to a foliation. Let $\; \mathcal U$ be a cocycle defining the foliation $\mathcal F$, $N_{\mathcal U}$ and   ${\mathcal H}_{\mathcal U}$ the associated transverse manifold and holonomy psudogroup, respectively.  First, by a transverse connection we understand a connection on the transverse manifold $N_{\mathcal U}$ of which elements of the holonomy pseudogroup ${\mathcal H}_{\mathcal U}$ are affine transformations. If such a connection exists for one transverse manifold, it exists on any other transverse manifold. A transverse connection $\hat{\nabla}$ defines a foliated connection $\nabla$  in the normal bundle by the formula

$$\bar{f}_i(\nabla_X\bar{Y} )= \hat{\nabla}_{\bar{f}_i(X)}\bar{f}_i(\bar{Y})$$

\noindent
for any vector $X \in TU_i$ and any foliated section $\bar{Y}$ on $U_i$. Conversely, by the same formula, any foliated connection defines  a connection on the transverse manifold of which the holonomy pseudogroup consists of affine transformations.

\medskip

Using the splitting of the tangent bundle $TM$ we can extend any connection in the normal bundle to a connection in $TM$ for which the subbundle $T{\mathcal F}$ is parallel. Conversely, any connection $' {\nabla}$  for which $T{\mathcal F}$ is parallel defines a connection $\nabla$ in the normal bundle by the formula below, where $\bar{..  }$ represents passing to the normal bundle (section):

$$\nabla_X\bar{Y} = \bar{'\nabla_XY}$$

\noindent
for any vector fields $X,Y.$ The induced connection $\nabla$ is foliated if for any infinitesimal automorphism (i.a.) of the foliation $X,$ and any  foliated section $\bar{Y}$ of the normal bundle, $\nabla_X\bar{Y}$ is a foliated setion of the normal bundle, which is equivalent to the fact that for any infinitesimal automorphism $Y$ of the foliation $\mathcal F$, $'\nabla_XY$ is an infinitesimal automorphism of $\mathcal F$.

\medskip
\noindent
{\bf Example 3} A foliation $\mathcal F$ is called transversely affine if there exsits a flat foliated connection. This is equivalent to the existence of a transverse flat connection. That is the transverse manifold is an affine (flat) manifold. Affine manifolds are locally affinely isomorphic to open subsets of $R^q$ with the standard flat connection. Therefore there exists a cocycle defining the foliation $\mathcal F$ modelled on the $R^q$ such that elements of the associated holonomy pseudogroup are restrictions of affine transformations of $R^q$. This is equivalent to the existence of an atlas adapted to the foliation such that the changes of the transverse coordinates are  restrictions of  affine transformations of $R^q$. Because that property transversely affine foliations are developable, cf. \cite{Wo1}.

\medskip

\noindent
{\bf Remark} Following \cite{Wo0}, we will use the  convention: "normal"  when qualifying a geometrical object means that this object is defined only on the normal bundle.  If such an object projects along local submersions defining the foliation, it will be called foliated. Holonomy invariant objects on a transverse manifold will be called transverse.  So foliated objects are normal but not all normal objects are foliated.  However, a connection $\nabla$ in the normal bundle $N(M,{\mathcal F})$ will be called normal if $\nabla_X\bar{Y} =0$ for any foliated section $\bar{Y}$ and $X \in T{\mathcal F}.$ This apparently inconsistency is only superficial   as only with that condition  a connection in the normal bundle can be defined as a section of a suitable associated bundle derived from the normal bundle.

 \section{Transversely Hessian foliations}

Let us continue with the study of a family of probability distributions as described in Introduction. Assume that  there exists a torsion-free connection $'\nabla ,$ cf. \cite{LL}, for which the tensor $g_F$ is parallel, i.e., $ '\nabla g_F = 0 ,$ and define  the  distribution $kerg_F :$

$$ kerg_F = \{ v \in TM \colon g_F(v,v) = 0\} = \{ v \in TM \colon g_F(v,w) = 0 \; \;  \forall w \in TM\}$$

\noindent
$kerg_F :$ is parallel with respect to $'\nabla $ as 

$$ 0='\nabla_X g_F(Y,Z) = Xg_F(Y,Z)  -  g_F('\nabla_XY,Z)  -  g_F(Y,'\nabla_XZ)$$

\noindent  
where $X,Y,Z \in TM.$  If $Y \in kerg_F,$ $g_F(Y,Z)=0$ for any $Z,$ and we get that $g_F('\nabla_XY,Z) = 0$ for any $Z \in TM,$ which ensures that $'\nabla_XY \in ker g_F$ provided that $Y \in kerg_F.$ As the connection $'\nabla $ is torsion-free the distribution $kerg_F$ is involutive:

$$0= T(X,Y) =' \nabla_XY -' \nabla_YX -[X,Y]$$

\noindent
so

$$[X,Y] = '\nabla_XY -' \nabla_YX \in kerg_F$$

\noindent
for any vector fields $X,Y \in kerg_F$. Moreover, the distribution $kerg_F$ is of constant dimension as $'\nabla_YX \in kerg_F$ for any $X \in kerg_F$ and any vector field $Y.$ This means that that the parallel transport along any curve maps vectors from $kerg_F$ to vectors from $kerg_F,$ which ensures that the distribution $kerg_F$ is of constant dimension, thus defines a foliation $\mathcal F$.

\medskip

 The tensor field $g_F$ induces a (normal) Riemannian metric $g$ in the normal bundle $N(M,{\mathcal F}).$ The connection $'\nabla$ defines a connection $\nabla$ in the normal bundle $N(M,{\mathcal F})$ which is the Levi-Civita connection of the metric $g.$  Indeed

$$\nabla_Xg(\bar{Y},\bar{Z}) = Xg(\bar{Y},\bar{Z})  - g(\nabla_X\bar{Y},\bar{Z}) 
- g(\bar{Y},\nabla_X\bar{Z}) = \;\;\;\;\;\;\;\;\;\;\;\;\;\;\;\;\;\;\;\;\;\;\;\;\;\;\;\;\;\;\;\;\;\;\;\;\;\;\;\;\;\;\;\;\;\;\;\;\;\;\;\;\;\;\;\;\;\;\;\;\;$$
$$ Xg_F(Y,Z)  - g(\bar{'\nabla_XY},\bar{Z}) 
- g(\bar{Y},\bar{'\nabla_X,Z}) = Xg_F(Y,Z)  - g_F('\nabla_XY,Z) - g_F(Y,'\nabla_X,Z) = 0$$

\noindent
where $X$ is any vector field, $Y,Z$ are i.a.s of the foliation. So the connection $\nabla$ is $g$-metric. Obviously, it is torsion-free, hence $\nabla$ is the Levi-Civita connection of the Riemannian metric $g$ of the normal bundle.


\medskip

The induced metric will be foliated if $Xg(\bar{Y},\bar{Z}) = Xg_F(Y,Z)=0$ for any vector field $X$ tangent to the foliation $\mathcal F$ and any foliated vector fields $Y$ and $Z$. It is the case if $'\nabla_XY \in T{\mathcal F}$  for any vector field $X$ tangent to the foliation $\mathcal F$ and any foliated vector fields $Y$.  As the connection $'\nabla$ is torsion-free this condition is equivalent to '$\nabla_YX \in T{\mathcal F}$ for any vector $Y\in TM$ and any vector field $X \in T{\mathcal F}.$

\medskip

The flat connection $'D$ of $\Lambda$ should be related in some way to the foliation $\mathcal F$. If we assume a similar condition for $'D,$  i.e. $'D_ZW \in T{\mathcal F}$ for any vector $Z\in TM$ and any vector field $W \in T{\mathcal F},$ then we can define a connection $D$ in the normal bundle by the formula

 $$D_X\bar{Y} = \bar{'D_XY}$$

\noindent
 where $\bar{Y}$ denotes the section of the normal bundle defined by the vector field $Y$. The normal connection $D$ is flat. The connection $D$ is transversally projectable (foliated) if $D_Xs$ is a foliated section for any foliated section $s$ and a foliated vector field $X,$ i.e. an infinitesimal automorphism of $\mathcal F$. 

\bigskip

Let $\mathcal F$ be a foliation of codimension $q$ on a manifold $M$ of dimension $m$. The dimension of its leaves is $p$, i.e. $p+q = m.$  Assume that the foliation $\mathcal F$ is transversely affine, cf. \cite{GHL,Wo1,Wo2}. A Riemannian metric $\hat{g}$ is bundle-like if  for any adapted chart $\varphi =(x^1,...,x^p,y^1,...,y^q),$  $\hat{g}=\Sigma_{ij=1}^pg_{ij}(x,y)v^iv^j + \Sigma_{a,b=1}^qg_{ab}(y)dy^ady^b,$ where $v^i$ is the only 1-form which vanishes on the orthogonal complement of the bundle $T\mathcal F$ and $v^i(\frac{\partial}{\partial x^j} =\delta^i_j.$  A bundle-like metric $\hat{g}$ is said to be transversely Hessian if  the horizontal part $g$ of the metric is expressed by the formula

$$g = \Sigma_{i.j=1}^{q} \frac{\partial^2 h}{\partial y^i \partial y^j}  dy^idy^j$$

\noindent
{\bf Remark}  In principle $h$ need not be basic, i.e. it may depend on variables $x^i$, we just assume that it does not.

\medskip

\noindent
{\bf Definition}
{\it 
We say that the foliation is transversely Hessian if 

- it is transversely affine

- it admits a bundle-like metric which is transversely Hessian.}

\medskip

\noindent
Therefore a transversely Hessian foliation is at the same time a Riemannian foliation and a transversely affine foliation althought both structures may have not much in common. 

\noindent
A foliation $\mathcal F$ is transversely Hessian iff for some cocycle  $\mathcal U$ defining  the foliation the associated transverse manifold $N_{\mathcal U}$ admits a Hessian structure of which elements of the holonomy pseudogroup are automorphisms, i.e. local affine transformations of the flat connection and isometries of the Hessian metric.

Let ${\nabla}$ be the (normal) Levi-Civita connection of the foliated Riemannian metric ${g}$ in $N(M, {\mathcal F})$. ${\nabla} $ is a foliated connection. The difference  ${\gamma} = {\nabla}  - {D}$ is also a normal tensor which is foliated iff ${D}$ is a foliated connection. Moreover, as both connections are torsion-free, ${\gamma}_XY = {\gamma}_YX.$ 

\medskip



Like in the standard case we have the following proposition

\begin{theorem} Let $(M,{\mathcal F})$ be a foliated manifold. Assume that ${\mathcal F}$ is transversely affine with foliated flat connection $D$ and $g$ is a foliated metric on $(M,{\mathcal F})$. Then, in any adapted chart $(U, \phi )$  the following conditions are equivalent, where $i,j,k$ run from 1 to q;

\medskip
i) $g$ is a foliated Hessian metric;

\medskip
ii) $D_Xg(Y,Z) = D_Yg(X,Z)$ for any foliated sections $X,Y,Z$ of the normal bundle;

\medskip

iii) $\frac{\partial g_{ij}}{\partial y^k} = \frac{\partial g_{kj}}{\partial y^i};$

\medskip

iv) $g(\gamma_XY,Z) = g(Y,\gamma_XZ)$  for any foliated sections $X,Y,Z$ of the normal bundle;

\medskip

v) $\gamma_{ijk} = \gamma_{jik}.$ 

\end{theorem}

The proof follows closely that Proposition 2.1 of \cite{Shima}. First of all, notice that as the metric $g$ foliated,  (i) implies (iii)  from the very definition of a Hessian metric. Let us choose and adapted atlas such that  the transverse coordinate changes are affine transformations. Then   the sections $\frac{\partial}{\partial y_i}$ are $D$-parallel. In such a coordinate system (iii) is just a local expression of (ii) and (v) of (iv). 

Moreover,  $\gamma_{ijk}$ are  related to the Christofel symbols of the foliated Levi-Civita connection $\nabla$. As such we have the following expression for them

$$\Gamma ^i_{jk} = ½ g^{is}( \frac{\partial g_{sj}}{\partial y_k} + \frac{\partial g_{sk}}{\partial y_j} -\frac{\partial g_{jk}}{\partial y_s})$$

\noindent
and 

$$\gamma _{ijk} =  g_{is}\Gamma ^s_{jk} =  ½ ( \frac{\partial g_{ij}}{\partial y_k} + \frac{\partial g_{ik}}{\partial y_j} -\frac{\partial g_{jk}}{\partial y_i})$$

\noindent
the derivatives with respect to variables along leaves of the foliation not appearing as the connection is foliated. It is clear that the conditions (iii) and (v) are equivalent. 

Let us demonstrate that (iii) implies (i). As the metric $g$ and the connection $D$ are foliated, and the local coordinate system is adapted to  the transversely affine foliation, the implication is equivatent to the corresponding fact for (transverse) manifolds. But it was proved in Proposition 2.1 of \cite{Shima}.

\medskip

Let $\mathcal F$ be a transversely affine foliation on a manifold M, let $N(M, {\mathcal F})$ be its normal bundle. It admits a foliation ${\mathcal F}_N$ of the same dimension as $\mathcal F$  but of codimension $2q.$

\medskip

Let $\varphi =(x^1,...,x^p,y^1,...,y^q)$ be an adapted chart. Then on the normal bundle $N(M, {\mathcal F})$ we have an adapted chart $\varphi^N = (x^1,...,x^p,y^1,...,y^q, dy^1,...,dy^q),$ 
of course the coordinates $x^i$ are identified with $x^ip$ where $p$ is the projection onto the base $M$ in the normal bundle. Moreover,  putting $z^i = y^i + \sqrt{-1} dy^i$ we define transversely holomorphic coordinate system on the normal bundle. Therefore if the foliations $\mathcal F$ is tranversely affine the foliation ${\mathcal F}_N$ is tranversely holomorphic for a foliated complex structure $J_N$. Moreover, on the normal bundle we define a normal Riemannian metric $g^N$ by the formula locally expressed

$$g^N = \Sigma_{ij}^q g_{ij}pdz^id\bar{z}^j$$

In this case we have the following proposition

\begin{proposition} Let $(M,{\mathcal F})$ be a foliated manifold, and $g$ be a foliated metric. Then the following conditions are equivalent:

(1) $g$ is a foliated Hessian metric on $(M,{\mathcal F},D)$

(2) $g^N$ is a foliated  K\"alerian metric on  $(N(M, {\mathcal F}), {\mathcal F}_N, J_N).$ 

\end{proposition}

From the very formulae both metrics are foliated. Both properties are local, so can be considered in a suitable adapted chart. In that case property (1) is equivalent to the fact that the induced transverse metric $g^T$ on $N_{\mathcal U}$  is Hessian, and property (2) to the fact that the induced metric ${g^N}^T$ on the tangent bundle of the manifold $N_{\mathcal U}$ is  K\"ahlerian. But it  is precisely the substance of Proposition 2.6 of \cite{Shima}.



\section{Dual foliated connections}


Let $(M,{\mathcal F}, D, g)$ be a tranversely Hessian foliated manifold. Let $\nabla$ be the Levi-Civita connection in the normal bundle of $\mathcal F$ associated to the foliated Riemannian metric $g$. Then the connection

$$ D' = 2\nabla  -  D$$

\noindent 
is a flat foliated  connection in the normal bundle and

$$Xg(Y,Z) = g(D_XY,Z) + g(Y,D'_XZ).$$

\noindent
Moreover, $(M,{\mathcal F}, D', g)$  is a transversely Hessian foliated manifold. $D'$ is called the dual connection of the  transversely Hessian foliated manifold $(M,{\mathcal F}, D, g).$


\medskip

Indeed, let us choose a cocycle $\; \mathcal U$  definfing the foliation. The associated transverse manifold $N_{\mathcal U}$ is a Hessian manifold of which the associated holonomy pseudogroup ${\mathcal H}_{\mathcal U}$ is a pseudogoup of local automorphisms.  Both connections $D$ and $\nabla$ are foliated, so they correspond to $\bar{D}$ and $\bar{\nabla},$  repectively. We can apply Corrolary 2.1 of \cite{Shima} to this Hessian structure with connections $\bar{\nabla}$ and $\bar{D}$ being the Levi-Civita connection and flat connection, respectively. The formula for the  dual connection $\bar{D}'$ (from Corollary 2.1  of \cite{Shima}) assures that any local automorphism of the Hessian structure is an affine transformation of the connection $\bar{D}'$, so it is a transverse connnection of our foliated structure. Therefore it defines a foliated connection $D'$ for which we have been looking for. 



\medskip

This situation can be well described using the notion of a Codazzi pair of connections. We say that a Riemannian metric $\bar{g}$ and a connection $\bar{D}$ in the normal bundle  are related by the normal Codazzi equations if

$${D}_{\bar{X}}{g}(\bar{Y},\bar{Z}) = {D}_{\bar{Y}}{g}(\bar{X},\bar{Z}) $$

\noindent
for any foliated sections $\bar{X},\bar{Y},\bar{Z}$ of the normal bundle.

\medskip

\noindent
{\bf Definition} {\it A  pair $(D,g),$  $g$ being a Riemannian metric in the normal bundle, and $D$ a torsion-free connection in this vector bundle, is called a normal Codazzi structure if it satisfies the 
Codazzi equation for  any i.a.s $X,Y$ and any foliated section $Z$:}

$$D_Xg(\bar{Y},Z) = D_Yg(\bar{X},Z).$$

\medskip
 If both objects are foliated, then the pair $(D,g)$ is called a {\it foliated Codazzi structure. } A  foliated manifold $(M,{\mathcal F})$ is called {\it a foliated Codazzi manifold,} and is denoted 
$(M,{\mathcal F};D,g).$ If $g$ is a foliated Riemannian metric and $D$ a foliated  connecction, then the pair is called 
a {\it  foliated Codazzi manifold. }

\medskip

Having a normal Codazzi structure  $({D},{g})$ on a foliated manifold $(M,{\mathcal F})$  we can define a new normal torsion-free connection ${D}'$ by the formula

$$X{g}(\bar{Y},\bar{Z}) = {g}({D}_X\bar{Y},\bar{Z}) + {g}(\bar{Y},{D}'_X\bar{Z})$$ 

\noindent
for any i.a.s $X,Y,Z$. The connection ${D}'$ is called the dual connection of the connection ${D}$ with respect to the Riemannian metric ${g}.$ 

\medskip

Let us check that the connection ${D}'$ is normal. If $X$ is tangent to the foliation then ${D}_X\bar{Y }=0$ as the connection ${D}$ is normal. Moreover, as $Y$ and $Z$ are i.a.s ${g}(\bar{Y},\bar{Z})$ is a basic function, so $X{g}(\bar{Y},\bar{Z}) =0.$ Thus ${g}(\bar{Y},{D}'_X\bar{Z})=0$ for any i.a. $Y, $  which ensures that in this case ${D}'_X\bar{Z}=0.$

\medskip

The connection $\bar{D}'$ is torsion-free as 

$${g}(\bar{Y},{D}'_X\bar{Z}- {D}'_Z\bar{X}) = {g}(\bar{Y},{D}_X\bar{Z}- {D}_Z\bar{X})$$

\noindent
for any i.a.s $X,Y,Z$. So ${D}'$ is torsion-free iff the connection ${D}$ is.

\medskip
Moreover, we have the following proposition

\medskip
\begin{proposition}  On a foliated manifold $(M,{\mathcal F}),$ if the Codazzi structure $(D,g)$  is foliated, so is the dual connection $D'. $
\end{proposition}

If the Codazzi structure $(D,g)$  is foliated, for any foliated sections $X,Y,Z$ the functions $Xg(Y,Z) $ and $g(D_XY,Z) $ are basic, so the function $g(Y,D'_XZ)$ is basic as well. That fact ensures that $D'_XZ$ is a foliated section for any foliated sections $X,Z$ of the normal bundle, which means that the connection $D'$ is foliated. 

\medskip

 In the case of normal connections we have the following lemma

\begin{lemma}
Let g be a foliated metric on the normal bundle. Let $(D,D')$ be a  pair of dual connections with respect to g.
Then if one of the connections is foliated so is the other.
\end{lemma}


The pair of  connections satisfies the equation

$$Xg(Y,Z) = g(D_XY,Z) + g(Y,D'_XZ).$$

\noindent
Let us check that $D'$ is a foliated connection if $D$ is. Let $Y,Z$ be two foliated section of $N(M,{\mathcal F}).$

\medskip

If $X$ is a vector field tangent to the leaves, then $Xg(Y,Z) =0$ and $D_XY=0$, so $D'_XY=0$ and the connection $D'$ is normal. 

\medskip

If $X$ is an i.a. of $\mathcal F$, then the function $Xg(Y,Z) $ is foliated (basic). Likewise the function $g(D_XY,Z)$ is foliated. Therefore $g(D'_XY,Z)$ is a foliated function, and the connection $D'$ is foliated. 

\medskip

 If $(M,{\mathcal F};D,g)$ is a foliated Codazzi manifold, so is $(M,{\mathcal F};D',g).$  In that case the induced connections $\bar{D}$ and $\bar{D}'$ on the transverse manifold are dual with respect to the transverse Riemannian metric $\bar{g}$ and the holonomy pseudogroup consists of affine transformations of both connections.

\medskip

The following lemma is a foliated version of a lemma for connections on manifolds, cf. Lemma 2.3 of  \cite{Shima}. The proof is basically the same. 

\begin{lemma}
Let $D$ be a torsion-free  connection and let $g$ be a Riemannian metric in the normal bundle. Define a new connection $D'$ by

$$Xg(Y,Z) = g(D_XY,Z) + g(Y,D'_XZ)$$

\noindent
for any vector field $X$ and any foliated sections $Y,Z$ of the normal bundle.  Then the following conditions (1)-(3) are equivalent:

(1) the connection $D'$ is torsion free,

(2)  the pair $(D,g)$ satisfies the Codazzi equation

$$D_Xg(Y,Z) = D_Yg(X,Z),$$

(3)  let $\nabla$ be the Levi-Civita connection for $g,$ and let $\gamma_XY = \nabla_XY-D_XY.$ Then $g$ and $\gamma$ satisfy

$$g(\gamma_XY,Z ) = g(Y,\gamma_XZ).$$

If the pair $(D,g)$ satisfies the Codazzi equation, then the pair $(D',g)$ also satisfies this equation and

$$D'=2\nabla -D,$$

$$D_Xg(Y,Z) = 2g ( \gamma_XY,Z)$$
\end{lemma}

\medskip
\noindent
{\bf Remark} Proposition 1 asserts that a flat connection $D$ and a Riemannian metric $g$ in the normal bundle of a foliation $\mathcal F$ form a  normal Hessian structure iff they satisfy the Codazzi equation. The same is true for foliated structures.

\medskip

The curvature tensor  $R_D$  of a connection $D$ in the normal bundle is defined as

$$R_D(X,Y)Z = D_XD_YZ -D_YD_XZ - D_{[X,Y]}Z$$

\noindent
for vector fields $X,Y$ and a section $Z$ of the normal bundle. If the connection $D$ is normal ($D_XZ =0$ for any vector tangent to leaves) the tensor field $R_D$ defines a tensor field denoted by the same letters $R_D$ a section of the bundle

$$N(M,{\mathcal F})^*\otimes N(M,{\mathcal F})^*\otimes N(M,{\mathcal F})^*\otimes  N(M,{\mathcal F})$$

\noindent
This bundle is also foliated, cf. \cite{Wo0}, and if the connection $D$ is foliated so the section $R_D$. If the connection $D$ is foliated the (foliated) curvature tensor field $R_D$ corresponds to  the curvature tensor field of the induced connection $\bar{D}$ on the transverse manifold.

\medskip

Therefore if the connection $D$ is foliated we can introduce the following definition

\begin{definition}A  foliated Codazzi structure is said to be of constant curvature c  if the curvature tensor $R_D$ of the connection $D$ satisfies

$$R_D(X,Y)Z = c ( g(Y,Z)X -g(X,Z)Y )$$

\noindent
for any sections $X,Y,Z$ of the normal bundle.

\end{definition}

\medskip

\noindent
{\bf Remark} The above condition is equivalent to the transverse Codazzi structure $(\bar{g}, \bar{D})$ being of constant curvature $c$ as  the curvature tensor $R_D$ of the connection $D$ projects, corresponds, to  the curvature tensor $R_{\bar{D}}$ of the connection $\bar{D}$ on the transverse manifold corresponding, induced by, the foliated connection $D$ . 

\medskip

The propositions below are  foliated versions of Propositions 2.8  and 2.9 of \cite{Shima}. They are immediate consequences of the just mentioned Propositions  applied to the induced Codazzi structures on the transvere manifold of the foliated manifold $(M,{\mathcal F})$ and of the above Remark.

\begin{proposition}

Let $(M,{\mathcal F}, D,g)$ be  a foliated Codazzi manifold, and $(D',g)$ be its dual Codazzi structure. 

(1) Denoting by $R_D$ and $R_{D'}$ the curvature tensors of $D$ and $D',$ respectively, we have

$$g(R_D(X,Y)Z,W) + g(Z,R_{D'}(X,Y)W) = 0,$$

\noindent
for any sections $X,Y,Z,W$ of the normal bundle.

(2) If $(D,g)$ is a Codazzi structure of constant curvature c, then $(D',g)$ is also a Codazzi structure of constant curvature c.

\end{proposition}

\begin{proposition}
A foliated Codazzi structure $(D,g)$ is of constant curvature 0 iff the foliation $\mathcal F$ is  transversely Hessian with the foliated structure given by the pair  $(D,g)$ .

\end{proposition}

\section{Normal curvatures}

Let $(D,g)$ be a normal Hessian structure on a foliated manifold $(M,{\mathcal F})$. Let $\gamma = \nabla -D$ be the difference tensor of the normal Levi-Civita connection of $g$ and $D$.

\medskip

A normal tensor field $Q$ of type (1,3) defined as 

$$Q= D\gamma$$

\noindent
is called the normal Hessian curvature tensor of  the normal Hessian structure $(D,g).$ 

\medskip

If the structure is foliated so is its normal Hessian curvature tensor. 
The components $Q^i_{jkl} $ of $Q$ with respect to an adapted foliated affine coordinate system
$(x^1,...,x^p,y^1,...,y^q) $ are given by, as the formula corresponds to the formula on the transverse manifold, see Chapter 3 of \cite{Shima}:

$$ Q^i_{jkl} = \frac{\partial \gamma ^i_{jl}}{\partial y^k}.$$

The propositions below are  foliated versions of Propositions 3.1 and 3.2 of \cite{Shima}:

\medskip

\begin{proposition}
Let us consider an adapted foliated affine coordinate system $(x^1,...,x^p,y^1,...,y^q) , $  and let $g_{ij} = \frac{\partial^2 \varphi }{\partial y^i y^j}$  for some basic function $\varphi .$ Then

(1) $ Q_{ijkl} =  \frac{1}{2}\frac{\partial^4 \varphi }{\partial y^i y^j y^k y^l }
- \frac{1}{2}g^{rs}\frac{\partial^3 \varphi }{\partial y^i y^k y^r}\frac{\partial^3 \varphi }{\partial y^j y^l y^s}.$

\medskip

(2) $Q_{ijkl} = Q_{klij}= Q_{kjil} = Q_{ilkj}= Q_{jilk}.$

\medskip

\noindent
for $i,j,k,l = 1,...,q.$
\end{proposition}

\begin{proposition}
Let $R$ be the Riemannian curvature tensor for the normal Riemannian metric $g.$ Then

$$ R_{ijkl} = \frac{1}{2} (Q_{ijkl} - Q_{jikl}).$$
\end{proposition}


\medskip

In Section  {\bf  Transversely Hessian Foliations}  we have demostrated that when $(M,{\mathcal F}) $ is transversely Hessian, the foliation ${\mathcal F}^N$ of the normal bundle $N(M,{\mathcal F}) $ is transversely K\"ahler for the lifted metric. The next proposition explains the relation between the tensor $Q$ and the curvature tensor of this K\"ahlerian metric.

\begin{proposition} Let $(M,{\mathcal F}) $ be transversely Hessian. Let $R^N$ be the Riemannian curvature tensor for the normal K\"ahlerian metric $g^N$ of the foliated manifold $(N(M,{\mathcal F}), {\mathcal F}^N).$ Then for any adapted foliated affine coordinate system

$$ R^N_{ijkl} = \frac{1}{2} Q_{ijkl}\circ \pi$$

\noindent
where $\pi$ is the natural projection $N(M,{\mathcal F}) \rightarrow M. $

\end{proposition}

It is a foliated vesion of Proposition 3.3 of \cite{Shima}, a simple consequence of the corresponce between "foliated" and "transverse"  objects and of the fact that the foliation ${\mathcal F}^N$ of the normal bundle  is defined by the cocycle derived from a cocycle defining the foliation $\mathcal F$ which is described at the beginning of section {\it Foliations}.

\medskip

A normal metric tensor defines a $q$-form $\omega$  by the standard formula

$$\omega_x(v_1,...,v_q) = det [g(v_i,v_j)]$$

\noindent
where $v_1,...,v_q \in N(M,{\mathcal F})_x.$  It can be understood as a section of $\wedge^q N(M,{\mathcal F})^*$.  The form $\omega$ is called  the normal volume form defined by $g$. If the metric $g$ is foliated, the the form $\omega$ is foliated (basic). The corresponding $q$ form $\bar{\omega}$  on the transverse manifold is the volume form defined by the transverse metric $\bar{g}$.  We define a closed 1-form $\alpha$ and a symmetric bilinear normal form $\beta$ by

$$ D_X\omega = \alpha (X) \omega,\;\;\;\; \beta = D\alpha .$$

\noindent
The forms $\alpha $ and $\beta$ are called the first normal Koszul form and the second   normal Koszul form, respectively, of the normal Hessian structure $(D,g)$ on a foliated manifold $(M,{\mathcal F})$. In the case of a foliated structure both forms are basic (foliated). The next proposition is a foliated version of Proposition 3.4  of \cite{Shima}. 

\begin{proposition} Let $(M,{\mathcal F}) $ be a  transversely Hessian foliated manifold.  Then for any adapted foliated affine coordinate system, with indices going from 1 to q,

(1) $\alpha (X) = Tr_N \gamma_X,$

(2) $\alpha_i = \frac{1}{2} \frac{\partial log det \vert g_{kl} \vert }{\partial y^i} = \gamma^r_{ri}$

(3)  $\beta_{ij} = \frac{\partial \alpha_i}{\partial y^j} = \frac{1}{2} \frac{\partial^2 log det \vert g_{kl} \vert  }{\partial x^i \partial y^j} = Q^r{_{rij}} =Q_{ij}{^r}{_r}.$
\end{proposition}

\medskip


The form $\beta$ is related to the normal  Ricci tensor ${\bf R}^N$ of the  normal K\"ahlerian metric $g^N$ of the foliated manifold $(N(M,{\mathcal F}), {\mathcal F}^N).$ The foliated version of Proposition 3.5 of \cite{Shima} reads as follows:

\begin{proposition}  Let $(M,{\mathcal F}) $ be transversely Hessian.  Then for any adapted foliated affine coordinate system $(x^1,...,x^p,y^1,...,y^q) $,  let ${\bf R}^N_{ij}$ be the local expression,with indeces going $ i,j$  from 1 to q, of  be the normal Ricci  tensor of the normal K\"ahlerian metric $g^N$ of the foliated manifold $(N(M,{\mathcal F}), {\mathcal F}^N).$ Then

${\bf R}^N_{ij} = - \frac{1}{2}\beta_{ij} \circ \pi .$

\end{proposition}


If the normal Hessian structure $(D,g)$ satisfies the condition

$$ \beta = \lambda g,  \;\;\; \lambda = \frac{{\beta^i}_i}{q}$$

\noindent
then the normal Hessian structure is called Einstein-Hessian. The theorem below explains the above condition and its relation to the well-known Einstein-K\"ahler condition for K\"ahlerian metrics. It is a direct consequence of the previous proposition.

\begin{theorem}
Let $(M,{\mathcal F})$ be a foliated manifold, and $(D,g)$ be a normal Hessian structure. Let $(J^N,g^N)$ be the normal K\"ahler structure of the  foliated manifold $(N(M,{\mathcal F}), {\mathcal F}^N)$ induced by $(D,g).$ Then the following conditions are equivalent

(1) $(M, {\mathcal F}, D,g) $ is   Einstein-Hessian;

(2) $(N(M,{\mathcal F}), {\mathcal F}^N, J^N, g^N)$ is  Einstein-K\"ahlerian.

\end{theorem}


\bigskip

\noindent
{\bf Normal Hessian sectional curvature}

\medskip

Let $Q$ be the normal Hessian curvature tensor. The formula, $i,j,k,l = 1,...,q,$ 

$$ \hat{Q}^{ ik} (\xi ) = \hat{Q}^{ i}{_j }{^k } {_l}\xi^{jl}$$

\noindent
defines an endomorphism $\hat{Q}$ on the space of symmetric contravariant normal two tensors, i.e., on $\otimes ^2N(M, {\mathcal F}).$ $\hat{Q}$ is self-dual (symmetric) with respect to the scalar product tensor induced by the Hessian metric $g$. 

\medskip

Let us define a function $q$ on the space $ \otimes ^2N(M, {\mathcal F})$ by the formula

$$q(\xi  ) = \frac{<\hat{Q}(\xi ), \xi > }{<\xi , \xi > } $$

\noindent 
for any $\xi \in \otimes ^2N(M, {\mathcal F})$ and $< , >$ the inner product defined by the normal Riemannian metric $g$ of the foliated Hessian structure. The function $q$ is called the normal Hessian sectional curvature. It is a basic (foliated) function if the normal Hessian structure is foliated. The manifold $ \otimes ^2N(M, {\mathcal F})$ admits the induced (natural) foliation of dimension $p$ whose leaves are coverings of leaves of the foliation $\mathcal F$, \cite{Wo0}. 

 We say that the foliated Hessian structure is of constant normal sectional curvature if $q$ is a constant function on $\otimes ^2N(M, {\mathcal F}).$ As the structures are foliated, and we have demonstrated that the discussed tensors are also foliated the results below, the  foliated versions of Proposition 3.6, Corollary 3.1 and 3.2  of \cite{Shima} are corollaries of these results when the corresponces between foliated and transverse geometric objects is applied.

\begin{proposition}
The normal Hessian sectional curvature of a foliated Hessian structure $(M,{\mathcal F},D,g )$ is constant and equal to c iff   for any adapted foliated affine coordinate system $(x^1,...,x^p,y^1,...,y^q) $,

$${Q}_{ijkl} = \frac{c}{2}\{{g}_{ij}{g}_{kl} + {g}_{il}{g}_{jk}\}$$
\noindent
for any $i,j,k,l = 1,...,q.$
\end{proposition}


\begin{corollary}Let $(M,{\mathcal F},D,g)$  be  a foliated Hessian structure.
Then the following two conditions are equivalent:

(1) the Hessian normal sectional curvature is a constant c;

(2) the holomorphic normal sectional curvature of $(N(M,{\mathcal F}), J^N,g^N;{\mathcal F}^N)$ is constant and equal to -c.

\end{corollary}


\begin{corollary}
Let $(M, {\mathcal F}, D,g)$  be  a foliated Hessian structure of constant normal Hessian sectional curvature c. Then the foliated Riemannian manifold $(M, {\mathcal F},g)$ is a foliated Riemannian manifold modeled on a space form of constant section curvature $-\frac{c}{4}.$
\end{corollary}







\section{Transversely Hessian foliations as developable foliations}

Let $X$ be a manifold and $G$ a subgroup  of $Diff(M)$ of diffeomorphisms of $X. $ We say that a foliation $\mathcal F$ of a smooth manifold $M$ is a $(X,G)$-foliation if it admits a cocycle $\mathcal U =\{ U_i , f_i , k_{ij} \}_{i\in I}$  modeled on $X$ such that the mappings $k_{ij}$ are restrictions of elements of the group $G$. Such a foliation is developable provided the action of $G$  is quasi-analytic, cf. \cite{GHL,Th},  i.e., if two diffeomorphisms of $G$ are equal on an open subset of $X,$ then they are equal. For example it is true for isometries of a Riemannian manifold. In such a case there exist

\medskip

i) a representation $h$ of the fundamental group of $M$ into $G$

$$ h \colon \pi_1(M, x_0) \rightarrow G$$

 ii) a developing mapping $D$ of  the universal covering $\tilde{M}$  onto $X$, i.e. a submersion 

$$D \colon \tilde{M} \rightarrow X$$

\noindent
which is $\pi_1(M, x_0)$- equivariant for the natural action on the universal covering space and the action on the manifold $X$ via the representation $h.$

iii) the fibres of $D$ define the foliation of $\tilde{M}$  which projects onto the foliation $\mathcal F$, i.e. the foliation by the connected componnents of the fibres of $D$ is the lift of the foliation $\mathcal F$ to  $\tilde{M}.$ 

Transversely Hessian foliations are developable as they are transversely affine, and as such they are $(R^q, Aff(R^q))$-foliations. Therefore the universal covering of $\Lambda$ admits a development $D$ onto $R^q$ with $imD$ being an open subset $R^q,$ which is invariant with respect to the action of $\pi_1(M, x_0)$ via the representation $h.$

Developable foliations have been studied in depth. It is easy to see that the transverse geometry can be read on $imD,$ and elements of $imh$ are automorphisms of these geometrical structures. So in the case transversely Hessian foliations on $imD$ in addition to the obvious flat connections we have a Riemannian metric of whose elements of $imh$ are isometries. In short, $imD$ is a Hessian manifold and elements of $imh$ are automorphisms of this Hessian structure.

\medskip

There are two salient questions to be answered:

a) is the developing mapping surjective?

b) are the fibres of $D$ diffeomorphic?

\medskip

The first question is a question about tranverse compleness of the foliation. There are no simple answers but for some see \cite{Wo2}. The second can be answered using the properties of the foliated Hessian metric. In fact.  the developing mapping $D$ is a Riemannian submersion, and  if the foliated Hessian metric is transversely complete, then $D$ is a locally trivial bundle. Local trivializations are obtained using geodesics orthogonal to the fibres. It is not easy to demonstrate transverse completeness of the metric. For example, it is the case when the foliated  manifold is compact.  Returning  to the initial example of a family of probability distributions parametrized by an open subset $\Lambda$ of $R^m$ it would be the case if $\Lambda$ admited a cocompact group $K$ of automorphisms of our transversely Hessian foliation, i.e. the quotient topological space ${\Lambda}/K$ is a compact manifold and the foliation $\mathcal F$ projects to thetransversely Hessian  foliation ${\mathcal F}_K$ of ${\Lambda}/K.$ Another possibility is the geodesic completeness of the foliated affine structure, cf. \cite{Wo2}, but even the xistence of a cocompact group of automorphisms does not ensure that this structure is geodesic complete. The results of the paper mentioned above point to the importance of the properties of the group $imh,$ called the holonomy group of the foliation. Its closure in $Aff(R^q$ and its Zariski closure play a particularly important role, see also \cite{GHL}. 

\medskip

The existence of two related transverse structures, Riemannian and affine, permits to combine the results of two well-developed theories. For both types of foliations we have good estimates of the growth type of leaves:

\medskip

- Riemannian foliations, see \cite{Ca,BG}

- transversely affine, see \cite{Wo1,Wo2,AW}

The growth type of leaves is bounded from above by the growth type of the holonomy group.

\medskip

The positive answer to the two questions formulated above definitely restricts the topology type of the manifold M, and the  domain $\Lambda$ of the principal example.

 The local triviality of the developing mapping $D$ ensures that the leaves of the foliation have diffeomorphic universal coverings provided that the manifold $M$ is connected.

\medskip

The existence of a transverse  invariant measure is assured by the fact that the foliation is Riemannian, for results and properties of invariant measures for transversely affine foliations see \cite{GHL}, also see \cite{Wo1}.


\begin{thebibliography}{99}
\bibitem{AW}  \'Alvarez L\'opez J. A. and  Wolak, R.: 
Growth of transversely homogeneous foliations. arXiv:1304.1939

\bibitem{AN} Amari S-I and Nagaoka H.:  Methods of Information Geometry,
Translations of Mathematical  monographs, AMS-OXFORD, vol 191 (2007)

\bibitem{BG} Breuillar, E. d and  Gelander T.: A topological Tits alternative, Ann. of Math. 166
(2007), 427 - 474.

\bibitem{CC} Candel, A. and Conlon L.: Foliations I,  Graduate Studies in Mathematics, vol. 23, AMS (2000)

\bibitem{Ca} Carri\`ere,  Y.:  Feuilletages riemanniens `a croissance polynomiale, Comment. Math.
Helv. 63 (1988), 1 - 20.

\bibitem{Fu} Fujimoto A.: Theory of G-strutures, Tokyo (1972)

\bibitem{GHL} Goldman, W., Hirsch, M. W., and Levitt, G.:  Invariant measures for affine foliations, 
Proc. AMS  86, 3, (1982),  511-518

\bibitem{Ha1}  Haefliger A.: Groupoïdes d'holonomie et classifiants, Astérisque, 116 (1984), 70 - 97

\bibitem{Ha2} Haefliger A.: Pseudogroups of local isometries. in: L.A. Cordero (Ed.), Differential Geometry, Santiago de Compostela, 1984, Res. Notes in Math., vol. 131, Pitman, Boston, MA (1985), 174 - 197

\bibitem{Ha3} Haefliger A.: Foliations and compactly generated pseudogroups. Foliations: Geometry and Dynamics, Warsaw, 2000, World Scientific, River Edge, NJ (2002), 275 - 295

\bibitem{LL}  Lehmann-Lejeune, J.: Int\'egrabilit\'e des $G$-structures d\'efinies par une 1-forme $0$-d\'eformable \`a  valeurs dans le fibr\'e tangent, Annales de l'Institut Fourier, 16,  no. 2 (1966), 329-387

\bibitem{Molino} Molino, P.: Riemannian Foliations.  Progress in Mathematics, no. 73,  Birkh\"auser
Boston, Inc., Boston, MA, (1988)

\bibitem{Shima} Shima, H.: The Geometry of Hessian  Structures. World Scientific (2007)

\bibitem{St} Sternberg, S.: Lectures on differential geometry. Prentice-Hall, Englewood Cliffs, N. J. (1964)

\bibitem{Th} Thurston, W.: The Geometry and Topology of Three-Manifolds. Princeton University Press (1997)

\bibitem{Wo0} Wolak, R.:  Foliated and associated geometric structures on foliated manifolds, Ann. Fac. Sci. Toulouse Math. (5) 10 (1989), no. 3, 337--360

\bibitem{Wo1} Wolak, R.: Transversely affine foliations compared with affine manifolds, Quart. J. Math. Oxford Ser. (2) 41 (1990), no. 163, 369--384

\bibitem{Wo2} Wolak, R.: Closure of leaves in transversely affine foliations, Canad. Math. Bull. 34 (1991), no. 4, 553--558

\bibitem{Wo3} Wolak, R.:
Growth of leaves in transversely affine foliations, Proc. Amer. Math. Soc. 127 (1999), no. 7, 2167--2173

\end{thebibliography}
\end{document}